\documentclass[12pt]{article}
\usepackage{fancyhdr}
\usepackage{amssymb}
\usepackage{amsmath}
\usepackage{amsthm}
\usepackage{graphicx}
\usepackage{psfrag}
\usepackage{graphics}
\usepackage{amsthm}

\input pictex.tex

\newcommand{\clo}{\mathrm{S}^1}

\theoremstyle{definition}

\newtheorem{thm}{Theorem}[section]
\newtheorem{prop}[thm]{Proposition}
\newtheorem{lem}[thm]{Lemma}
\newtheorem{rem}[thm]{Remark}

\newtheorem{ex}[thm]{Example}

\begin{document}

\pagestyle{fancy}
\rhead{On conjugates and distortion of diffeomorphisms}

\date{}
\author{Andr\'es Navas}

\title{On conjugates and the asymptotic distortion of 1-dimensional $C^{1+bv}$ diffeomorphisms}

\maketitle

\vspace{-0.1cm}

\noindent{\bf Abstract:} We show that a $C^{1+bv}$ circle diffeomorphism with absolutely continuous derivative 
and irrational rotation number can be conjugated to diffeomorphisms that are $C^{1+bv}$ arbitrary close to the 
corresponding rotation. This improves a theorem of M.~Herman, who established the same result but starting 
with a $C^2$ diffeomorphism. We prove that the same holds for countable Abelian groups of circle diffeomorphisms 
acting freely, a result that is new even in the $C^{\infty}$ context. Related results and examples concerning the 
asymptotic distortion of diffeomorphisms are presented. Along this path, we provide a straightened 
version of the classical Denjoy-Kocsma inequality for absolutely continuous potentials.

\vspace{0.2cm}


\vspace{0.2cm}


\vspace{0.75cm}

\hfill {\em In memory of Jean-Christophe Yoccoz}

\vspace{0.25cm}

\noindent{\Large{\bf Introduction}}

\vspace{0.35cm}

Since the seminal work of Poincar\'e and Denjoy, circle diffeomorphisms became a central object of study in dynamics. 
One of the most celebrated results of the theory is the linearization theorem of Herman \cite{herman}, which 
establishes the differentiability of the conjugating map to the corresponding rotation in the case of an irrational 
rotation number satisfying a Diophantine condition. The key idea of Herman's proof consists in conjugating the 
diffeomorphism into one that is arbitrarily close to the rotation in a good topology, and finally applying local 
linearization results coming from standard KAM techniques (see also \cite{herman-simple}). 

Herman's theorem was later extended via different approaches by many people, including 
Yoccoz \cite{yoccoz}, Katznelson and Orstein \cite{KO}, and Khanin and Teplinsky \cite{KT}. 
However, the problem of approximation by conjugation raised by Herman in Chapter VII of his seminal 
work \cite{herman} is interesting by itself. In this regard, just a couple of results that apply to diffeomorphisms of arbitrary 
irrational rotation number are known. On the one hand, it was shown by the author in \cite{conj} that every $C^1$ circle diffeomorphism 
$f$ of irrational rotation number $\rho = \rho (f)$ has a sequence of conjugates that converges in $C^1$ topology to the rotation $R_{\rho}$; 
actually, there is a continuous path of $C^1$ conjugates that starts at $f$ and finishes at $R_{\rho}$. On the other hand, Herman proved in 
\cite{herman} that if $f$ is of class $C^2$, then there exists a sequence of conjugates of $f$ that converges to $R_{\rho}$ in the $C^{1+bv}$ 
sense;  that is, they converge to $f$ in the uniform topology, their derivatives uniformly converge to 1, and the total variation of the logarithm 
of these derivatives converge to 0. The main goal of this paper is to relax the regularity hypothesis of $f$ for this result and, simultaneously, 
extend it to commuting diffeomorphisms. 

\vspace{0.5cm}

\noindent{\bf Main Theorem.} {\em Let $\, \Gamma$ be a countable Abelian group of $C^{1+bv}$ (orientation preserving) 
circle diffeomorphisms whose action is free. If all elements in $\Gamma$ have absolute continuous derivative, then there exists a path 
$g_t$, $t \in [0,1[$, of $C^{1+bv}$ diffeomorphisms starting at $g_0 = id$ such that, for each $f \in \Gamma$, the conjugates $\, g_t f g^{-1}_t \,$ 
converge to the rotation $R_{\rho (f)}$ in the $C^{1+bv}$ sense as $t$ goes to 1.}

\vspace{0.5cm}

A sequence of good conjugates for the statement above can be made explicit for $\mathbb{Z}^k$-actions. 
Since these will be frequently used, we state a general proposition on this.

\vspace{0.5cm}

\noindent{\bf Main Proposition.} 
{\em In the statement of the Main Theorem above, assume that $\Gamma$ is isomorphic to $\mathbb{Z}^k$, and 
let $\{ f_1, \ldots, f_k\}$ be a generating system. For each $n$, denote $B(n)$ the set of elements 
of the form $f_1^{n_1} f_2^{n_2} \cdots f_k^{n_k}$ for which $0 \leq n_i \leq n$, and 
let $h_n$ be a diffeomorphism such that 
\begin{equation}\label{explicit-general}
D h_n (x) = 
\frac{\left[ \prod_{g  \in B(n-1)} Dg  (x) \right]^{1 / |B(n-1)|}}{\int_{\mathrm{S}^1} \left[\prod_{g \in B(n-1)} Dg (z) \right]^{1 / |B(n-1)|} d z}.
\end{equation}
Then, for each $f \in \Gamma$, the conjugates  
$h_n f h_n^{-1}$ converge to $R_{\rho(f)}$ in the $C^{1+bv}$ sense.}

\vspace{0.5cm}

There are several aspects involved in these two statements, as we next discuss.

\vspace{0.25cm}

\noindent {\bf I.} First notice that though the Main Theorem is stated for countable groups, it easily reduces to finitely-generated ones, 
as we will explain at the end of the proof. The reason for stating it in such a generality relies on that it applies to some relevant groups, 
as for instance those constructed by Hirsh in \cite{Hirsh} and exploited in \cite{GS} and \cite{Na-kova}.

\vspace{0.25cm}

\noindent {\bf II.} Assume now that $\Gamma$ is isomorphic to $\mathbb{Z}$, hence generated by a single (infinite order) 
circle diffeomorphism $f$ with absolutely continuous derivative. If the action of $\Gamma$ is free, then 
$f$ has irrational rotation number, and our Main Theorem yields 
the existence of a path of conjugates $g_tfg_t^{-1}$ that converges to the rotation $R_{\rho}(f)$ in the $C^{1+bv}$ sense.  
Actually, according to the Main Proposition, a sequence of such conjugates can be made explicit: letting $h_n$ be such that 
\begin{equation}\label{explicit-Z}
Dh_n (x) = \frac{\left[ \prod_{k=0}^{n-1} Df^k (x) \right]^{1/n}}{\int_0^1 \left[ \prod_{k=0}^{n-1} Df^k (z) \right]^{1/n} dz},
\end{equation}
we have the desired convergence for $h_n f h_n^{-1}$. 

There is a direct consequence of this fact that relies on a deep theorem obtained by Yoccoz in his 
PhD thesis (that solved in the affirmative a conjecture of Herman \cite{herman}): linearizable 
circle diffeomorphisms of irrational rotation number are dense (in the $C^{\infty}$ topology). 

\vspace{0.5cm}

\noindent{\bf Corollary 1.} Conjugacy classes of $C^{1+bv}$ circle diffeomorphisms of irrational rotation number and absolute 
continuous derivative are $C^{1+bv}$ dense. 

\vspace{0.18cm}

\noindent{\bf Proof.} Let us start by making more precise the statement. It claims that, given any two $C^{1+bv}$ circle diffeomorphisms $f,g$ 
with the same irrational rotation number and having absolutely continuous derivative, there exists a sequence of $C^{1+bv}$ diffeomorphisms 
$h_n$ such that the conjugates $h_n f h_n^{-1}$ converge to $g$ in the $C^{1+bv}$ sense. This means that $h_n f h_n^{-1}$ converges to $g$ 
in the $C^1$ topology and the total variation of the logarithm of the derivative of both $g^{-1} (h_n f h_n^{-1})$ and $(h_n f h_n^{-1}) g^{-1}$ go 
to zero as $n$ goes to infinity.

To show this, first notice that $g$ is $C^{1+bv}$ close to a $C^{\infty}$ diffeomorphism $\hat{g}$ of the same rotation number 
$\rho$. (This is easy to check for a $C^2$ diffeomorphism, but it is slightly tricker -~and it is left as an exercise- for a $C^{1+bv}$ diffeomorphism  
with absolute continuous derivative.) By Yoccoz' theorem (see \cite[Chapitre 3]{yoccoz-ast}), $\hat{g}$ is $C^{\infty}$ close to a diffeomorphism 
$\tilde{g}$ that is $C^{\infty}$ linearizable, say $\tilde{g} = h R_{\rho}h^{-1}$ for a $C^{\infty}$ diffeomorphism $h$. By the Main Theorem, there 
exists a sequence of $C^{1+bv}$ diffeomorphisms $h_n$ such that $h_n f h_n^{-1}$ converges to $R_{\rho}$ in the $C^{1+bv}$ sense. It is 
then easy to see that $h (h_nfh_n^{-1}) h^{-1}$ converges to $h R_{\rho} h^{-1} = \tilde{g}$ in the $C^{1+bv}$ sense. Since $\tilde{g}$ 
can be taken arbitrarily close to $\hat{g}$ and this last map can be taken arbitrarily close to $g$, 
this gives the desired sequence of conjugates of $f$ converging to $g$. $\hfill\square$

\vspace{0.5cm}

\noindent{\bf III.} A third aspect to comment on the Main Theorem concerns the proof itself. Let us first 
mention that the strategy is similar to that of \cite{conj}. Namely, we look for a translation of the  problem of 
approximation by conjugates into a problem of finding approximate solutions to a certain cohomological 
equation. In the $C^1$ case, the later concerns the logarithmic 
derivative and is a classical cohomological equation. In the present setting, the equation that arises deals with 
the affine derivative and is twisted by the action of the classical derivative. The absolute continuity 
assumption translates into that this affine derivative takes values in $L^1$, which makes 
it very pleasant to treat by using functional analytical methods.

At this point, one would wonder whether a similar approach could be implemented for the space of diffeomorphisms having derivatives in 
$L^p$ for some $p > 1$, which is nicely suited in other contexts (as for instance that of linearization problems; see \cite{KO}). Unfortunately, 
in our setting, this would involve computations that are beyond the scope the methods employed here, hence no result in this direction is 
provided here. Notice that the case $p=\infty$ corresponds to diffeomorphisms with Lipschitz derivative, which is very close to the $C^2$ 
case for which our result discussed in V. requires a  very strong hypothesis.

Next we pass to the key ingredient of the proof. This is nothing but a standard Hahn-Banach type argument, which  
shows that the existence of approximate solutions for the 
twisted cohomological equation in $L^1$ follows from the ergodicity with respect to Lebesgue measure. This is nicely complemented 
by the fact  that ergodicity of $C^{1+bv}$ circle diffeomorphisms of irrational rotation number is a classical result of Katok \cite{KH},  
also proved by Herman for $C^{1}$ diffeomorphisms with Lipschitz derivative \cite{herman}. 

There is a quite unexpected consequence of this observation that involves Birkhoff sums of absolutely continuous functions.  
This may be seen as a complement to the Denjoy-Kocsma inequality. It extends a result from \cite{NT}, which was established 
for $C^1$ functions using  the uniqueness of the 1-conformal measure proved by Douady and Yoccoz in \cite{DY}. 

\vspace{0.5cm}

\noindent{\bf Theorem 1.}
{\em Let $f$ be a $C^{1+bv}$ circle diffeomorphism of irrational rotation number.  If $\phi$ is an absolutely continuous function 
on the circle then, in $C^0$ topology, one has the convergence 
$$\lim_{n \to \infty} \left[ \phi + \phi \circ f + \ldots + \phi \circ f^{q_n - 1} - q_n \int_{\clo} \phi \, d \mu \right] = 0,$$ 
where $q_n$ denotes the sequence of denominators of the rational approximations of $\rho (f)$ and $\mu$ denotes the 
unique invariant probability measure of $f$.}

\vspace{0.5cm}

Recall that a map $f$ is said to be {\em rigid} if there exists an increasing sequence of integers $q_n$ such that $f^{q_n}$ converges to the identity 
(in a prescribed topology). A fundamental result of Herman \cite{herman} establishes that every $C^2$ circle diffeomorphism $f$ of irrational rotation 
number is $C^1$ rigid, with $q_n$ being the sequence of denominators for the rational approximations of the rotation number of $f$. By applying 
Theorem 1 above to the function $\phi := \log (Df)$, we obtain the next slight extension of this result.

\vspace{0.5cm}

\noindent{\bf Corollary 2.} {\em Every $C^{1+bv}$ circle diffeomorphism $f$ of irrational rotation number that has absolutely continuous derivative
is $C^1$ rigid. More precisely, the sequence $f^{q_n}$ converges to the identity in $C^1$ topology, where $q_n$ is the sequence of denominators 
for the rational approximations of the rotation number of $f$.}

\vspace{0.5cm}

\noindent{\bf Proof of Corollary 2.} First recall that $\log (Df)$ has zero mean with respect to the unique $f$-invariant probability measure 
(see \cite{herman}). Theorem 1 then implies that 
$$\log (Df^{q_n}) = \sum_{i=0}^{q_n - 1} \log (Df) (f^{i})$$
uniformly converges to zero as $n$ goes to infinity. Therefore, $Df^{q_n}$ uniformly converges to 1. 
Since $f^{q_n}$ converges to the identity in $C^0$ topology (because of Denjoy's theorem), this proves the $C^1$ convergence 
of $f^{q_n}$ to the identity. $\hfill\square$

 \vspace{0.5cm}

\noindent{\bf IV.} A fourth relevant aspect of the Main Theorem relies to another byproduct of the proof. Namely, a straightforward argument 
shows that the existence of almost solutions to the twisted cohomological equation associated to the affine derivative of a circle 
diffeomorphism $f$ implies that the affine derivatives of the powers $f^n$ grow sublinearly. This means that the value of 
 $$\frac{\mathrm{var}(\log(Df^n))}{n}$$ 
 converges to $0$ as $n$ goes to infinite. We then consider this notion for general diffeomorphisms of one-dimensional manifolds. 
 More precisely, we call the {\em asymptotic distortion} of such an $f$ the number 
$$\mathrm{dist}_{\infty} (f) := \lim_{n \to \infty} \frac{\mathrm{var}(\log(Df^n))}{n}.$$ 
For the case of  circle diffeomorphisms of irrational rotation number, the panorama summarizes  in the following result. 

\vspace{0.5cm}

\noindent{\bf Theorem 2.} 
{\em The asymptotic distortion of every circle diffeomorphism of irrational rotation number and absolutely continuous derivative 
always vanishes. However, there exist $C^{1+bv}$ circle diffeomorphisms of irrational rotation number whose asymptotic distortion 
does not vanish.}

\vspace{0.5cm}

The next result clarifies the relation between asymptotic distortion and closures of conjugacy classes of diffeomorphisms.

\vspace{0.5cm}

\noindent{\bf Theorem 3.} 
{\em  If $f$ is a $C^{1+bv}$ diffeomorphism of a compact one-dimensional manifold $M$, then its asymptotic distortion 
vanishes if and only if there exists a sequence $h_n$ of $C^{1+bv}$ diffeomorphisms such that the conjugates 
$h_n f h_n^{-1}$ converge to an isometry in the $C^{1+bv}$ sense.}
 
 \vspace{0.5cm}
 
Some other examples related to asymptotic distortion will be exhibited. In particular, we will prove 
that, in $\mathrm{PSL}(2,\mathbb{R})$, only the hyperbolic elements have positive asymptotic distortion. 
However, we will show that there exist  $C^{\infty}$ diffeomorphisms of the interval with positive asymptotic distortion 
but without hyperbolic fixed points. Actually, we will see that this is the typical behavior of diffeomorphisms with no hyperbolic fixed points. 

Let us mention that the notion of asymptotic distortion should be extendable to diffeomorphisms of higher 
dimensional manifolds. Although the behavior seems both worth and interesting to study in this setting, 
it is much less clear. A key difficulty consists in the failure of invariance of  the variation of the logarithm 
of the derivative under change of coordinates. (This is actually the case of any other reasonable funcional \cite{bruno}.)

\vspace{0.5cm}

\noindent{\bf V.} A last aspect to discuss on the Main Theorem involves higher regularity. Indeed, 
although part of our approach still applies in the $C^2$ setting, a key technical point (concerning 
the ``size'' of the conjugating map) at the end of the proof fails, so that the question of whether 
conjugates of $C^2$ circle diffeomorphisms of arbitrary irrational rotation number can approximate 
the corresponding rotation in $C^2$ topology remains open. (This is Question 25 of \cite{Rio}; 
it should be pointed out that Avila and Krikorian have recently announced a result in this direction 
but for $C^{\infty}$ conjugates of $C^{\infty}$ diffeomorphisms.) The theorem below points in this direction; 
it extends a result from Chapter VII of \cite{herman}, which holds for single diffeomorphisms. We stress, 
however, that the hypothesis that we assume (namely, $C^1$ linearizability) is a very strong one.

\vspace{0.5cm}

\noindent{\bf Theorem 4.} {\em If $\Gamma$ is a finitely-generated Abelian group of $C^2$ circle diffeomorphisms that is conjugate 
to a group of rotations by a $C^1$ diffeomorphism, then there exists a sequence $h_n$ of $C^2$ diffeomorphisms such that, for all 
$f \in \Gamma$, the sequence of conjugates $g_nfg_n^{-1}$ converges to a rotation in $C^2$ topology as $n$ goes to infinity.}

\vspace{0.5cm}

The plan of this paper is as follows.  
In \S \ref{section-key}, we start by proving a Key Lemma on approximate solutions to the twisted cohomological equation associated to the affine 
derivative of a circle diffeomorphism. We next proceed to establish some direct consequences of this: Theorem~1 is proved in \S \ref{section-DK},  
the vanishing of the asymptotic distortion for $C^1$ diffeomorphisms with absolutely continuous derivative is proved in \S \ref{section-VAD}, 
and the existence of $C^{1+bv}$  diffeomorphisms with nonvanishing asymptotic distortion is proved in \S \ref{section-NVAD} (all these 
results concern diffeomorphisms of irrational rotation number, and the last two ones putted together yield Theorem 2).  In \S \ref{section-main}, 
we prove the Main Theorem and the Main Proposition.  In \S \ref{section-AD}, we come back to the notion of asymptotic distortion; in particular, 
we prove Theorem  3.  Finally, \S \ref{section-C2} is devoted to the proof of Theorem 4. 
 
We stress that all maps considered in this article will be assumed to be orientation preserving.

\section{A key lemma on twisted cohomological equations} 
\label{section-key}

We immediately proceed to the Key Lemma 
of this article, which is somewhat related to (though simpler than) the main result of \cite{NT}. After this, we pass to the 
results announced in the Introduction that can be directly deduced from this.

\vspace{0.5cm}

\noindent{\bf Key Lemma.} {\em Let $f$ be a $C^{1+bv}$ circle diffeomorphism of irrational rotation number. If $\phi$ is 
an $L^1$ function defined on the circle that has zero mean w.r.t Lebesque measure, then there exists a sequence of $L^1$ 
functions $\psi_n$ of zero mean such that the sequence}
$$\psi_n - \psi_n \circ f \cdot Df$$
{\em converges to $\phi$ in the $L^1$ sense as $n$ goes to infinity.} 

\vspace{0.2cm}

\noindent{\bf Proof.} Assume first that $\phi$ cannot be approximated by functions of the form 
\, $\psi - \psi \circ f \cdot Df$ \, in $L^1$,  
where no condition on the mean of $\psi \in L^1$ is imposed. 
Then, by the Hahn-Banach theorem, there exists a linear functional $L \!: L^1 (\clo) \to \mathbb{R}$ 
that vanishes on the set of functions of the form $\psi - \psi \circ f \cdot D f$ and 
such that $L (\phi) = 1$. This functional comes from 
integration against an $L^{\infty}$ function $\xi$:
$$L (\zeta) = \int_{\clo}  \xi \cdot \zeta.$$
The equality $L (\psi - \psi \circ f \cdot D f) = 0$ then becomes
$$\int_{\clo} \xi \cdot \psi = \int_{\clo} \xi \cdot \psi \circ f \cdot Df = \int_{\clo} \xi \circ f^{-1} \cdot \psi,$$
where the last equality follows from a change of variable. 
Since this holds for every $\psi \in L^1$, this implies that a.e. we have the equality $\xi = \xi \circ f^{-1}$. By the 
Katok-Herman ergodicity theorem, the function $\xi$ is a.e. constant, say equal to $c$. However, this implies that
$$L(\phi) = c \int_{\clo} \phi = 0,$$
which is a contradiction. 

Let hence $\tilde{\psi_n}$ be a sequence of $L^1$ functions such that the sequence
$$\tilde{\psi}_n - \tilde{\psi}_n \circ f \cdot Df$$
converges to $\phi$ in $L^1$. For each $n$ let 
\, $c_n := \int \tilde{\psi}_n.$ \, 
Then the function $\tilde{\psi}_n - c_n$ has zero mean, and satisfies 
$$(\tilde{\psi}_n - c_n) - (\tilde{\psi}_n - c_n) \circ f \cdot Df + c_n (1 - Df) = \tilde{\psi}_n - \tilde{\psi}_n \circ f \cdot Df \longrightarrow \phi.$$
To deal with the unpleasant expression $c_n (1-Df)$, let 
$$\zeta_k := \frac{1}{q_k} \left[ 1 + Df + Df^2  + \ldots + Df^{q_k-1} \right] - 1,$$
where $q_k$ is the sequence of denominators in the rational approximation of $\rho(f)$. 
Notice that $\zeta_k$ has zero mean for all $k \geq 1$. Moreover, an easy computation yields 
$$1 - Df = (\zeta_k - \zeta_k \circ f \cdot Df) - \frac{1}{q_k} \big[ Df^{q_k} - 1 \big].$$
By Denjoy's inequality \cite{libro}, the last expression on the right uniformly converges to zero, hence 
\, $\zeta_k - \zeta_k \circ f \cdot Df$ \,
uniformly converges to $1 - Df$.  The desired sequence $\psi_n$ is then obtained by letting 
$$\psi_n := \tilde{\psi}_n - c_n + c_n \zeta_{k_n}$$ 
for a well-chosen sequence $k_n$. 
$\hfill\square$


\subsection{An application to the Denjoy-Kocsma inequality}
\label{section-DK}

The goal here is to prove Theorem 1, that we recall for the reader's convenience.

\vspace{0.5cm}

\noindent{\bf Theorem 1.}
{\em Let $f$ be a $C^{1+bv}$ circle diffeomorphism of irrational rotation number.  If $\phi$ is an absolutely continuous function 
on the circle then, in $C^0$ topology, one has the convergence 
$$\lim_{n \to \infty} \left[ \phi + \phi \circ f + \ldots + \phi \circ f^{q_n - 1} - q_n \int_{\clo} \phi \, d \mu \right] = 0,$$ 
where $q_n$ denotes the sequence of denominators of the rational approximations of $\rho (f)$ and $\mu$ denotes the unique 
invariant probability measure of $f$.}

\vspace{0.5cm}

\noindent{\bf Proof.} 
Since $\phi$ is absolutely continuous, its derivative $D \phi$ is an $L^1$ function, and it has zero mean w.r.t. Lebesgue measure. 
By the Key Lemma, given $\varepsilon > 0$, there exists function $\psi \in L^1$ with zero mean w.r.t. Lebesgue measure such that the $L^1$ 
function $\xi$ defined by
$$\xi := D\phi - \big( \psi - ( \psi \circ f) \cdot Df \big)$$
has $L^1$-norm at most $\varepsilon/2$. Notice that $\xi$ also has zero mean w.r.t. Lebesgue measure. By integrating, we can view $\psi$ and 
$\xi$ as the derivatives of functions $\Psi$ and $\Xi$ on the circle, respectively. The equality above then reads as 
$$D \Xi = D \phi - (D \Psi - D (\Psi \circ f)).$$
Notice that $\Xi$ and $\Psi$ are defined up to a constant. We fix the representative of $\Xi$ with zero mean w.r.t. $\mu$. Then 
the equation above becomes 
\begin{equation}\label{zero-mean-mu}
\Xi = \phi - (\Psi - \Psi \circ f) - \int_{\clo} \phi \, d\mu.
\end{equation}
A telescopic sum then yields
\begin{equation}\label{telescopic}
\sum_{i=0}^{q_n - 1} \phi \circ f^i - q_n \int_{\clo} \phi \, d\mu = \sum_{i=0}^{q_n-1} \Xi \circ f^i + (\Psi - \Psi \circ f^{q_n}).
\end{equation}
Now, by the Denjoy-Kocsma inequality \cite{herman,libro}, 
$$\left\| \sum_{i=0}^{q_n-1} \Xi \circ f^i \right\|_{C^0} \leq \mathrm{var} (\Xi) = \| D \Xi \|_{L^1} = \| \xi \|_{L^1} \leq \frac{\varepsilon}{2}.$$
Moreover, since $\Psi$ is a continuous function and $f^{q_n}$ converges to the identity in $C^0$ topology (because of Denjoy's theorem), 
for a large-enough $n$ we have \, $\| \Psi - \Psi \circ f^{q_n} \|_{C^0} \leq \varepsilon / 2$. \, By (\ref{telescopic}), a triangular inequality 
yields, for a large $n$,  
$$\left\| \sum_{i=0}^{q_n - 1} \phi \circ f^i  - q_n \int_{\clo} \phi \, d\mu \right\|_{C^0} \leq \varepsilon.$$
Since this holds for every $\varepsilon > 0$, this shows the announced convergence. $\hfill\square$


 \subsection{Vanishing of the asymptotic distortion}
 \label{section-VAD}
 
 As we announced, the existence of approximate solutions to the twisted cohomological equation implies the vanishing of the 
 asymptotic distortion. This corresponds to the first half of Theorem 2.
 
\begin{prop} \label{prop-ADV}
If $f$ is a $C^{1+bv}$ circle diffeomorphism with absolutely continuous detivative and 
irrational rotation number, then the sequence
$$\frac{1}{n} \left\| \frac{D^2 f^n}{D f^n} \right\|_{L^1} = \frac{\mathrm{var} (\log (Df^n); \clo)}{n}$$ 
converges to 0 as $n$ goes to infinity.
\end{prop}

\vspace{0.1cm}

\noindent{\bf Proof.} Denote by $U \!: L^1 \to L^1 $ the linear isometry given by $U(\varphi) = \varphi \circ f \cdot Df$. 
The cocycle relation
\begin{equation}\label{cociclo}
\frac{D^2 (g_1 g_2)}{D (g_1 g_2)} =  \frac{D^2 g_2}{D g_2} + \frac{D^2 g_1}{D g_1} \circ g_2\cdot D g_2 
\end{equation}
yields 
\begin{equation}\label{ene}
\frac{D^2 f^n}{D f^n} = \sum_{i=0}^{n-1} U^i \left( \frac{D^2 f}{D f} \right).
\end{equation}
Recall that $D^2 f / Df$ is an $L^1$ function with zero mean w.r.t. Lebesgue measure. 
By the Key Lemma, given $\varepsilon > 0$, there exists an $L^1$ function $\psi$ such that 
$$\left\| \frac{D^2 f}{D f} - (\psi - \psi \circ f \cdot Df) \right\|_{L^1} \leq \varepsilon.$$
Since $U$ is an isometry, this implies that, for all $i \geq 0$, 
\begin{small}
$$\left\| U^i \left( \frac{D^2 f}{D f} \right) - (\psi \circ f^i \cdot Df^i - \psi \circ f^{i+1} \cdot Df^{i+1}) \right\|_{L^1} = 
\left\| U^i \left( \frac{D^2 f}{D f} - (\psi - \psi \circ f \cdot Df) \right) \right\|_{L^1} \leq \varepsilon.$$
\end{small}Using the triangular inequality together with (\ref{ene}), this gives 
\begin{small}
$$\left\| \frac{D^2 f^n}{Df^n} - (\psi - \psi \circ f^n \cdot Df^n) \right\|_{L^1} 
\leq \sum_{i=0}^{n-1} \left\| U^i \left( \frac{D^2 f}{D f} \right) - (\psi \circ f^i \cdot Df^i - \psi \circ f^{i+1} \cdot Df^{i+1}) \right\|_{L^1} 
\leq n \varepsilon.$$
\end{small}Therefore,
$$\frac{1}{n} \left\| \frac{D^2 f^n}{D f^n} \right\|_{L^1} 
\leq \varepsilon + \frac{1}{n} \Big( \|\psi\|_{L^1} + \| (\psi \circ f^n) \cdot Df^n \|_{L^1} \Big) 
= \varepsilon + \frac{2 \| \psi \|_{L^1} }{n}.$$
Since this holds for every $\varepsilon > 0$, the announced convergence easily follows.
$\hfill\square$


\subsection{On low regular diffeomorphisms}
\label{section-NVAD}

Our goal here is to prove the second half of Theorem 2.

\vspace{0.15cm}

\begin{prop}
\label{prop-NVAD}
There exist $C^{1+bv}$ circle diffeomorphisms of irrational rotation number with nonvanishing asymptotic distortion.
\end{prop}

\vspace{0.15cm}

Before passing to the proof, recall that Mather constructed in \cite{mather} a surjective, continuous group homomorphism 
$\mathcal{M}$ on the group of $C^{1+bv}$ circle diffeomorphisms that takes values in $\mathbb{R}$. 
His construction proceeds as follows: given $f \!\in\! \mathrm{Diff}^{1+bv}_+ (\clo)$, let $\mu_f$ be the signed measure 
obtained as the representative (under the Riesz representation theorem) of the Riemann-Stieljes integration with respect to $\log (D f)$. 
By the Hahn decomposition theorem, this decomposes as the sum of an absolute continuous part $\mu_f^{ac}$ and a totally singular 
part $\mu_f^{s}$. The value of Mather's homomorphism $\mathcal{M} \!: \mathrm{Diff}^{1+bv}_+(\clo) \to \mathbb{R}$ at $f$ is then defined as 
$$\mathcal{M} (f) := \int_{\clo} d \mu_f^{s} = - \int_{\clo} d \mu_f^{ac}.$$
Notice that an $f \in \mathrm{Diff}^{1+bv}_+ (\clo)$ can be easily constructed with any prescribed rotation number and any prescribed 
value for $\mathcal{M}(f)$.
 
\vspace{0.25cm}

\noindent{\bf Proof of Proposition 1.2.} If $f \!\in\! \mathrm{Diff}^{1+bv}_+ (\clo)$ satisfies $\mathcal{M} (f) \neq 0$, then it cannot hold that  
\begin{equation}\label{conv}
\frac{1}{n} \mathrm{var} (\log (Df^n)) \longrightarrow 0
\end{equation}
as $n$ goes to infinity. Indeed, one has $\mathcal{M} (f^n) / n = \mathcal{M} (f)$ for all $n$, and the continuity of $\mathcal{M}$ 
recalled above directly implies that, if (\ref{conv}) holds, then $\mathcal{M} (f^n) / n$ converges to zero. $\hfill\square$

\vspace{0.15cm}

\begin{rem}\label{rem-mina}
Notice that Mather's homomorphism can be extended to piecewise-differentiable maps for which the logarithm of the 
derivative has finite total variation. However, in the framework of piecewise-affine homeomorphisms, 
his construction does not apply. More 
precisely, the measure $\mu_f$ induced by the logarithm of the derivative is totally singular, hence Mather's homomorphism is trivial. 
Nevertheless, the asymptotic distortion may still fail to vanish, and this failure can be completely detected by 
a single invariant for the case of maps with irrational rotation number. 
Namely, letting $Df_+$ (resp. $Df_{-}$) be the right (resp. left) derivative and denoting $J_f (x) := Df_+ (x) / Df_{-} (x)$ 
the {\em jump} of $f$ at $x$, the crucial ingredient in this description is the complete jump along orbits:
$$ C_f (x) := \prod_{n \in \mathbb{Z}} J_f (f^n(x)).$$
Notice that this product is well defined, since it involves only finitely many factors that are different from $1$. 

On the one hand, according to Minakawa \cite{minakawa}, if $C_f (x)$ equals 1 at every point, 
then $f$ is conjugate to the corresponding rotation by either a piecewise-affine homeomorphism or by such a map 
together with a homeomorphism of the form $x \mapsto \frac{\sigma^x - 1}{\sigma - 1}$, which is piecewise real-analytic with a derivative of finite 
total variation. In each case, the sequence $\mathrm{var} (\log (Df^n); \clo)$ is bounded. 
On the other hand, if there exists $x$ for which $C_f (x) \neq 1$, 
then choosing cut points $x_i$ very close (from the left and from the right) to points of the form $f^{-k} (x)$ (for $k$ large enough so that 
$J_f (f^{-k-\ell}(x)) = 1$ for all $\ell \geq 0$), one easily checks that, along such a partition of the circle, the variation of the logarithm of 
the derivative of $f^n$ is larger than  \, $n |\log (C_f (x))|$ \, minus a constant. Therefore, 
$$\lim_{n \to \infty} \frac{\mathrm{var} (\log (Df^n); \clo)}{n} \geq | \log (C_f (x_0)) | > 0.$$ 
\end{rem}

\vspace{0.25cm}

We close this section with a couple of remarks and questions. First notice that, 
following the lines of \cite[Chapter IV]{herman} (see also \cite{np}), 
it is not hard to prove that, for a $C^{1+bv}$ circle diffeomorphism (independently 
of the rotation number), the sequence $\mathrm{var} (\log (Df^n); \clo)$ remains bounded if and only if $f$ is conjugated to the corresponding 
rotation by a $C^{1+bv}$ diffeomorphism. A slightly subtler remark is that if, moreover, $f$ has 
absolutely continuous derivative, then the conjugacy can be taken with absolutely continuous derivative as well. (This follows as an 
application of the $L^1$ fixed point theorem from \cite{BGM}, and holds for any group action on the circle along which $D^2 g / D g$ 
remains uniformly bounded in $L^1$.) The next question is inspired by the growth gap established in \cite{watanabe}.

\vspace{0.5cm}

\noindent{\bf Question 1.} In case of unboundedness, does there exist a constraint for the growth of the sequence 
$\mathrm{var} (\log (Df^n); \clo)$ stronger than sublinearity that depends on the regularity of $f$?

\vspace{0.5cm}

The next question is inspired by Theorem 1 and the proof of Proposition \ref{prop-NVAD}.

\vspace{0.5cm}

\noindent{\bf Question 2.} Does there exist a $C^{1+bv}$ circle diffeomorphism with irrational rotation number that is not $C^1$ rigid? 
Can such a diffeomorphism be detected using Mather's homomorphism?


\section{Proof of the Main Theorem}
\label{section-main}

The Main Theorem can be established for $\mathbb{Z}$-actions as a consequence of the Key Lemma. However, the case of general Abelian 
groups of diffeomorphisms requires several adjustments. This is the reason why we start with a conceptual lemma that somewhat reduces the 
general case to that of $\mathbb{Z}$-actions. This is stated in a broad conceptual framework for potential use in different contexts. For 
our purposes, it will be used for the isometry $U \!: \varphi \mapsto ( \varphi \circ f) \cdot Df$ of $L^1$ and the associated cocycle
$$c(f) = \frac{D^2 f}{D f}.$$

\vspace{0.015cm}

\begin{lem} \label{lema-U}
Let $U$ be a linear isometric action of a finitely-generated Abelian group $\Gamma$ on a 
Banach space $\mathbb{B}$, and let $c \!: \Gamma \to \mathbb{B}$ be a cocycle, that is, a function that satisfies the relation 
\begin{equation}\label{coc-id}
c (g_1 g_2) = c(g_2) + U(g_2) (c(g_1)).
\end{equation} 
Assume that, for every element $g$ of infinite order in $\Gamma$, one has 
$$\frac{\| c(g^n) \|_{\mathbb{B}}}{n} \rightarrow 0$$
as $n \to \infty$. Then there exists a sequence of  functions $\psi_n \in \mathbb{B}$ such that, for 
all $f \in \Gamma$, the sequence $\psi_n - U(f) (\psi_n)$ converges to $c(f)$ as $n$ goes to infinite.
\end{lem}

\vspace{0.2cm}

\noindent{\bf Proof.} Let us first assume that $\Gamma$ is finitely generated and has no torsion, say $\Gamma \sim \mathbb{Z}^k$. 
Let  $\{f_1,\ldots, f_k\}$ be a (minimal) system of generators of $\Gamma$. For each $n \in \mathbb{N}$, 
let us denote by $B(n)$ be the {\em positive square ball} of radius $n$, that is, 
\, $B (n) := \{ f_{1}^{n_1} f_2^{n_2} \cdots f_k^{n_{k}}: \,\, 0 \leq n_i \leq n \}.$ \, 
For each $n \geq 1$, define
\begin{equation}
\psi_n := \frac{1}{|B(n-1)|} \sum_{g \in B(n-1)} c(g).
\label{int-zero}
\end{equation}
Then, for each $f \in \Gamma$, 
\begin{eqnarray*}
U(f) (\psi_n) 
&=& \frac{1}{| B(n-1) |} \sum_{g \in B(n-1)} U(f) ( c(g) ) \\
&=&  \frac{1}{| B(n-1) |} \sum_{g \in B(n-1)} [c(gf) - c(f)] \\
&=& - c(f) + \frac{1}{ | B(n-1) | } \sum_{g \in B(n-1)} c(g f) \\
&=& - c(f) + \frac{1}{ | B(n-1) | } \sum_{g \in B(n-1)} c(f g),
\end{eqnarray*}
where we have used the cocycle identity (\ref{coc-id})  in the second equality above. 
Therefore, for each $i \in \{1,2,\ldots,k\}$,
$$\Big\| c ( f_i ) - ( \psi_n - U(f_i) (\psi_n) ) \Big\|_{\mathbb{B}} \leq \frac{1}{B(n-1)} \Big\| \sum_{g \in B(n-1)} [ c (f_i g) - c(g) ] \Big\|_{\mathbb{B}},$$
and the last expression equals
$$\frac{1}{B(n-1)} \left\| \sum_{\substack{0 \leq n_j \leq n-1 \\ j \neq i}} \big[ c (f_i^n f_1^{n_1} \cdots f_{i-1}^{n_{i-1}} f_{i+1}^{n_{i+1}} \cdots f_k^{n_k}) 
- c (f_1^{n_1} \cdots f_{i-1}^{n_{i-1}} f_{i+1}^{n_{i+1}} \cdots f_k^{n_k}) \big] \right\|_{\mathbb{B}}.$$
By the cocycle relation, this reduces to
$$\frac{1}{B(n-1)} \left\| \sum_{\substack{0 \leq n_j \leq n-1 \\ j \neq i}} 
U (f_1^{n_1} \cdots f_{i-1}^{n_{i-1}} f_{i+1}^{n_{i+1}} \cdots f_k^{n_k}) (c (f_i^n) ) \right\|_{\mathbb{B}},$$
which, by the triangular inequality, is smaller than or equal to 
$$\frac{1}{B(n-1)} \sum_{\substack{0 \leq n_j \leq n-1 \\ j \neq i}} 
\left\| U (f_1^{n_1} \cdots f_{i-1}^{n_{i-1}} f_{i+1}^{n_{i+1}} \cdots f_k^{n_k}) (c (f_i^n) ) \right\|_{\mathbb{B}} 
= \frac{1}{B(n-1)} \sum_{\substack{0 \leq n_j \leq n-1 \\ j \neq i}} \| c (f_i^n) \|_{\mathbb{B}}.$$
A simple counting argument yields the following upper bound for the last expression:
$$\frac{1}{n^k} \, n^{k-1} \| c (f_i^n) \|_{\mathbb{B}} = \frac{\| c (f_i^n) \|_{\mathbb{B}}}{n},$$
According to the hypothesis, this goes to zero as $n \to \infty$. Therefore,
$$\| c ( f_i ) - ( \psi_n - U(f_i) (\psi_n) ) \|_{\mathbb{B}} $$
converges to zero as $n \to \infty$. To show that this holds changing $f_i$ by any $f \in \Gamma$, just write $f$ as 
a product of generators, apply the previous conclusion to each of these generators, and use the triangular inequality. 

Assume now that $\Gamma$ is still finitely generated but has torsion, 
say $\Gamma \sim \mathbb{Z}^k \times T$, where $T$ is a finite group. Let $\psi_n^*$ be 
a sequence for which $\| c ( f ) - ( \psi_n^* - U(f) (\psi_n^*) ) \|_{\mathbb{B}} $ converges to zero for every $f \in \mathbb{Z}^k$. Replace 
$\psi_n^*$ by 
\begin{equation}
\psi_n := \frac{1}{|T|} \sum_{g \in T} \, [ U (g) (\psi_n^*) + c(g) ].
\label{int-zero2}
\end{equation}
For each $f \in \Gamma$, we have
\begin{eqnarray*}
U(f) \psi_n \!
&=& \frac{1}{|T|} \sum_{g \in T} [ U(f) U(g) (\psi_n^*) + U(f) c(g) ] \\
&=& \frac{1}{|T|} \sum_{g \in T} [ U(f) U(g) (\psi_n^*) + c (gf) - c (f)] 
\,\,\, = \,\,\,-c(f) +  \frac{1}{|T|} \sum_{g \in T} [ U(gf) (\psi_n^*) + c(gf) ]. 
\end{eqnarray*}
If $f$ belongs to $T$, then the last expression reduces to $-c(f) + \psi_n$, thus showing that we actually have an equality:
$$c(f) = \psi_n - U(f) (\psi_n).$$
Assume otherwise that $f$ belongs to $\mathbb{Z}^k$. Then 
\begin{small}
\begin{eqnarray*}
U(f) \psi_n 
= -c(f) +  \frac{1}{|T|} \sum_{g \in T} [ U(gf) (\psi_n^*) + c(gf) ] 
= -c(f) +  \psi_n + \frac{1}{|T|} \sum_{g \in T} [ U(gf) (\psi_n^*) + c(gf) ] - \psi_n,
\end{eqnarray*}
\end{small}hence the value of 
$$\| c(f) - (\psi_n - U(f) \psi_n) \|_{\mathbb{B}}$$
is bounded from above by
\begin{small}
$$\frac{1}{|T|} \Big\|  \sum_{g \in T} 
[ U(gf) (\psi_n^*) + c(gf) ] - [ U(g) (\psi_n^*) + c(g) ] \Big\|_{\mathbb{B}} 
= \frac{1}{|T|} \Big\|  \sum_{g \in T}
U(g) [U(f) (\psi_n^*) - \psi_n^*] + [c(fg) - c(g)] \Big\|_{\mathbb{B}},$$
\end{small}that is, by
\begin{eqnarray*}
\frac{1}{|T|} \Big\|  \sum_{g \in T} U(g) [U(f) (\psi_n^*) - \psi_n^*] + U(g) (c(f)) \Big\|_{\mathbb{B}} 
&=& 
\frac{1}{|T|} \Big\|  \sum_{g \in T} U(g) [U(f) (\psi_n^*) - \psi_n^* + c(f) ] \Big\|_{\mathbb{B}} \\
&=&
\Big\|  U(f) (\psi_n^*) - \psi_n^* + c(f) \Big\|_{\mathbb{B}},
\end{eqnarray*}
which, by construction, converges to zero as $n \to \infty$. 
Now, to deal with an arbitrary $f \in \Gamma$, just write $f = g_1 g_2$, with $g_1 \in \mathbb{Z}^k$ and $g_2 \in T$, 
and use the cocycle relation and the triangular inequality to show the desired convergence for the function $\psi_n$. 

Finally, let us consider the case of a countable yet non-necessarily  finitely-generated group $\Gamma$. Let 
$\mathcal{G}_1 \subset \mathcal{G}_2 \subset \ldots$ be a complete axhaustion of $\Gamma$ by finite subsets and, for each 
$n$, let $\Gamma_n$ be the group generated by $\mathcal{G}_n$. We have already established the claim of the Lemma for each 
$\Gamma_n$. As a consequence, for each $n$, there exists a function $\psi_n$ such that, for all $f \in \mathcal{G}_n$, one has 
$$\big\|c (f) - (\psi_n - U(f) (\psi_n)) \big\|_{\mathbb{B}} < \frac{1}{n}.$$
Since $\bigcup_n \mathcal{G}_n = \Gamma$, this sequence $\psi_n$ clearly fulfils the desired property. 
$\hfill\square$

\vspace{0.15cm}

\begin{rem} 
The proof above mimics one of the key arguments of \cite{conj}.  (Compare \cite{Bo-Na}.) 
Actually, some of the ideas therein allow extending it to (locally) nilpotent group actions.
\end{rem}

\vspace{0.15cm}

We are now in position to prove our Main Theorem.

\vspace{0.5cm}

\noindent{\bf Proof of the Main Theorem.} Let us consider the isometries $U \!: \varphi \mapsto ( \varphi \circ f) \cdot Df$ 
defined on the space $L^1_0$ made of the $L^1$ functions with zero mean, as well as the associated cocycle
$$c(f) = \frac{D^2 f}{D f} \in L^1_0.$$
(Notice that (\ref{coc-id}) corresponds to (\ref{cociclo}) in this framework.) Since the action is free, infinite-order elements 
have irrational rotation number. Proposition \ref{prop-ADV} then establishes the hypothesis of Lemma \ref{lema-U}, which 
yields a sequence of functions $\psi_n \in L^1_0$ such that, for all $f \in \Gamma$, the value of  
$$\| c(f) - (\psi_n - \psi_n \circ f \cdot Df )\|_{L^1}$$
converges to zero as $n \to \infty$. Let $h_n$ be defined by 
\begin{equation}\label{def-psi-ene}
h_n (x) := \frac{ \int_0^x \exp \Big( \int_0^y \psi_n (z) \, dz \Big) \, dy }{ \int_0^1 \exp \Big( \int_0^y \psi_n (z) \, dz \Big) \, dy }.
\end{equation}
Since $\psi_n$ has zero mean, this map $h_n$ is a $C^1$ circle diffeomorphism. (The condition on the mean 
translates into that $Dh_n (0) = Dh_n (1)$.) Besides, in the $L^1$ sense, we have the equality
$$\frac{D^2 h_n}{D h_n} = D (\log (Dh_n)) = \psi_n.$$
For each $f \in \Gamma$, the cocycle relation (\ref{cociclo}) implies
$$\frac{D^2 (h_n f h_n^{-1})}{D (h_n f h_n^{-1})} 
= \frac{D^2 h_n}{D h_n} \circ (f h_n^{-1} ) \cdot D (f h_n^{-1}) 
+ \frac{D^2 f}{D f} \circ h_n^{-1} \cdot Dh_n^{-1}
- \frac{D^2 h_n}{D h_n} \circ h_n^{-1} \cdot D h_n^{-1}.$$ 
Therefore, 
$$\frac{D^2 (h_n f h_n^{-1})}{D (h_n f h_n^{-1})} 
= \left[ \psi_n \circ f \cdot D f  + \frac{D^2 f}{D f} - \psi_n \right] \circ h_n^{-1} \cdot D h_n^{-1}.$$ 
Besides, given $\varepsilon > 0$, there exists $N \in \mathbb{N}$ such that, for all $n \geq N$, we have 
$$\left\| \psi_n \circ f \cdot D f  + \frac{D^2 f}{D f} - \psi_n \right\|_{L^1} \leq \varepsilon.$$
As a consequence, by change of variable,
$$\left\| \Big[ \psi_n \circ f \cdot D f  + \frac{D^2 f}{D f} - \psi_n \Big] \circ h_n^{-1} \cdot D h_n^{-1} \right\|_{L^1} \leq \varepsilon.$$ 
Therefore,
$$\mathrm{var} (\log (D(h_n f h_n^{-1}))) 
= \left\| \frac{D^2 (h_n f h_n^{-1})}{D (h_n f h_n^{-1})} \right\|_{L^1} 
\leq \varepsilon,$$
which implies that $h_n f h_n^{-1}$ converges to the rotation of angle $\rho(f)$ in the $C^{1+bv}$ sense. 

Finally, to build the desired path of conjugating maps $g_t$, just interpolate the functions $\psi_n$ above. More 
precisely, let $u: [0,1[ \to [0,\infty[$ be a homeomorphism. For $t \in [0,1[$, let 
$n$ be the integer such that $n \leq u(t) < n+1$, and define
\, $\psi^t := ( n+ 1  - u(t) ) \, \psi_n + ( u(t) - n) \, \psi_{n+1},$\ \,
where $\psi_0 \equiv 0$. Finally, define the diffeomorphism $g_t$ by letting
$$g_t (x) := \frac{ \int_0^x \exp \Big( \int_0^y \psi^t (z) \, dz \Big) \, dy }{ \int_0^1 \exp \Big( \int_0^y \psi^t (z) \, dz \Big) \, dy }.$$
We leave to the reader the task of checking that the sequence of conjugates $g_t f g_t^{-1}$ satisfies the desired properties. 
$\hfill\square$

\vspace{0.3cm}

Following the lines of the proofs above, the conjugating maps can be made explicit. The case of $\mathbb{Z}^k$-actions is 
particularly pleasant and is summarized in our Main Proposition recalled bellow. For the general case, one would just 
need to ``make the conjugators homogeneous along torsion elements", in the spirit of the proof of Lemma \ref{lema-U}.

\vspace{0.5cm}

\noindent{\bf Main Proposition.} 
{\em In the framework of the Main Theorem, assume that $\Gamma$ is isomorphic to $\mathbb{Z}^k$, and 
let $\{ f_1, \ldots, f_k\}$ be a generating system. For each $n$, denote $B(n)$ the set of elements 
of the form $f_1^{n_1} f_2^{n_2} \cdots f_k^{n_k}$ for which $0 \leq n_i \leq n$, and let $h_n$ be 
a diffeomorphism such that  
$$
D h_n (x) = 
\frac{\left[ \prod_{g \in B(n-1)} Dg  (x) \right]^{1 / |B(n-1)|}}{\int_{\mathrm{S}^1} \left[\prod_{g \in B(n-1)} Dg (z) \right]^{1 / |B(n-1)|} d z}.
$$
Then, for each $f \in \Gamma$, the conjugates  
$h_n f h_n^{-1}$ converge to $R_{\rho(f)}$ in the $C^{1+bv}$ sense.}

\vspace{0.5cm}

\noindent{\bf Proof.} In this setting, definition (\ref{int-zero}) becomes 
$$\psi_n = \frac{1}{|B(n-1)|} \sum_{g \in B(n-1)} \frac{D^2 g}{D g}.$$
A simple computation then shows that, for $h_n$ given by (\ref{def-psi-ene}), the value of $Dh_n (x)$ equals
$$\frac{1}{\int_{0}^{1} \! \left[ \prod_{g \in B(n-1)} \frac{Dg (z)}{Dg (0)} \right]^{ 1 / B(n-1)} \! dz} 
\left( \prod_{g \in B(n-1)} \frac{Dg (x)}{Dg (0)} \right)^{\!\! 1 / B(n-1)} 
\! = \, 
\frac{\left[ \prod_{g \in B(n-1)} Dg (x) \right]^{1/B(n-1)}}{\int_0^1 \! \left[ \prod_{g \in B(n-1)} Dg (z) \right]^{1/B(n-1)} dz}.$$
To conclude, notice that this computation only uses that $D^2 h_n / D h_n = \psi_n$. Therefore, the same equality holds if we change $h_n$ 
by a left composition of it with a rotation, and this corresponds to just prescribing the value of its derivative as in (\ref{explicit-general}). $\hfill\square$

\vspace{0.5cm}

The reader will notice that the maps $h_n$ coincide with the conjugating maps considered in \cite{conj}. For the case of a 
$\mathbb{Z}$-action generated by a diffeomorphism $f$, they satisfy
$$Dh_n (x) = \frac{\left[ \prod_{k=0}^{n-1} Df^k (x) \right]^{1/n}}{\int_0^1 \left[ \prod_{k=0}^{n-1} Df^k (z) \right]^{1/n} dz}.$$
Remark that these are quite different from the conjugating maps used by Herman given by
$$h_n^* (x) := \frac{1}{n} \sum_{k=0}^{n-1} ( f^k - k \rho(f)),$$
where, on the right, $f$ denotes a lift to the real-line of the original map, and $\rho(f)$ is its corresponding translation number. 
It is worth mentioning that certain analogs of Herman's conjugating maps $h_n^*$ still work for Abelian (and even nilpotent) 
group actions in the continuous framework, yet providing explicit estimates for the derivatives of the conjugate action in the 
smooth case seems out of reach. (Herman managed to control derivatives of $\mathbb{Z}$-actions induced by $C^2$ 
diffeomorphisms.)


\section{On the asymptotic distortion of diffeomorphisms}
\label{section-AD}

Given a $C^{1+bv}$ diffeomorphism of a compact one-dimensional manifold $M$, let us recall that its {\em asymptotic distortion} is defined as  
$$\mathrm{dist}_{\infty} (f) := \lim_{n \to \infty} \frac{\mathrm{var} (\log (Df^n); M)}{n} = \lim_{n \to \infty} \frac{\mathrm{var} (\log (Df^n))}{n},$$
where, on the right and henceforth, we drop the reference to the manifold $M$ to simplify the notation. 
Notice that the existence of the limit follows from subaditivity:
$$\mathrm{var} (\log (Df^{m+n})) \leq \mathrm{var} (\log (Df^m)) + \mathrm{var} (\log (Df^n)).$$
In a certain sense, $\mathrm{dist}_{\infty}$ is a ``stable version'' of the variation of the logarithm of the derivative, 
in the sense that it satisfies, for all $k \geq 1$,
$$\mathrm{dist}_{\infty} (f^k) = k \cdot \mathrm{dist}_{\infty} (f).$$
Moreover, it is a dynamical invariant. Indeeed, the subaditivity relation 
$$\mathrm{var} (\log(D(fg))) \leq \mathrm{var} (\log(D f)) + \mathrm{var} (\log(D g))$$ 
easily implies that $\mathrm{dist}_{\infty}$ is invariant under $C^{1+bv}$ conjugacy, since
$$ \mathrm{var} (\log (Df^k)) + 2  \mathrm{var} (\log (Dh)) \leq 
\mathrm{var} (\log( D( (hfh^{-1})^k))) \leq \mathrm{var} (\log (Df^k); M) + 2 \mathrm{var} (\log (Dh)).$$

Theorem 3 (whose statement is recalled below) somewhat summarizes the relation between 
the vanishing of the asymptotic distortion and the approximation (in the $C^{1+bv}$ sense) of an 
isometry (either a rotation for the circle case or the identity for the interval case) by a sequence of conjugates. 
Notice that, in this theorem, no assumption on absolute continuity for the derivative is made. 

\vspace{0.5cm}

\noindent{\bf Theorem 3.} 
{\em  If $f$ is a $C^{1+bv}$ diffeomorphism of a compact one-dimensional manifold $M$, then its asymptotic distortion 
vanishes if and only if there exists a sequence $h_n$ of $C^{1+bv}$ diffeomorphisms such that the conjugates 
$h_n f h_n^{-1}$ converge to an isometry in the $C^{1+bv}$ sense.}

\vspace{0.2cm}

\noindent{\bf Proof.} If $h_n f h_n^{-1}$ converges to an isometry then, for each $\varepsilon > 0$, there exists an integer $N$ 
such that, for all $n \geq N$, 
$$\mathrm{var} (\log (D (h_n f h_n^{-1}) )) < \varepsilon.$$
This yields, for every integer $k$,
$$\mathrm{var} (\log(D (h_n f^k h_n^{-1}))) 
= \mathrm{var} (\log (D ((h_n f h_n^{-1})^k))) 
\leq k  \, \mathrm{var} (\log (D(h_n f h_n^{-1})))  
< k \, \varepsilon.$$
Therefore, 
$$\mathrm{var} (\log(Df^k)) 
\leq
\mathrm{var} (\log(D(h_n f^k h_n^{-1}))) + 2 \, \mathrm{var} (\log (D h_n)) 
< 
k \, \varepsilon + 2 \, \mathrm{var} (\log (Dh_n)),$$
hence
$$\frac{\mathrm{var} (\log(Df^k))}{k} < \varepsilon + \frac{2 \, \mathrm{var} (\log (Dh_n))}{k}.$$
Letting $k \to \infty$, we obtain
$$\lim_{k \to \infty} \frac{\mathrm{var} (\log(Df^k))}{k} \leq \varepsilon.$$
Since this holds for all $\varepsilon > 0$, the asymptotic distortion of $f$ vanishes.

The converse statement cannot be proved following the lines of the proof of the Main Theorem, since this would require absolutely 
continuous derivatives. Instead of pursuing this functional analytical approach, we use the explicit conjugating maps $h_n$ given 
by (\ref{explicit-Z}). Letting $y:= h_n^{-1}(x)$, a straightforward computation using the chain rule yields
\begin{eqnarray*}
D ( h_n f h_n^{-1} ) (x)
&=&
\frac{D h_n (f(y)) }{D h_n (y)} \cdot Df (y)
\,\,\, = \,\,\,
\frac{\left[ \prod_{k=0}^{n-1} Df^k (f(y)) \right]^{1/n}}{\left[ \prod_{k=0}^{n-1} Df^k (y) \right]^{1/n}} \cdot Df (y) \\
&=&
\frac{\left[ \prod_{k=0}^{n-1} Df^k (f(y)) \cdot Df (y) \right]^{1/n}}{\left[ \prod_{k=0}^{n-1} Df^k (y) \right]^{1/n}} 
\,\,\, = \,\,\,
\frac{\left[ \prod_{k=1}^{n} Df^k (y) \right]^{1/n}}{\left[ \prod_{k=0}^{n-1} Df^k (y) \right]^{1/n}} 
\,\,\, = \,\,\,
Df^n (y)^{1/n}.
\end{eqnarray*}
Therefore, 
$$\mathrm{var} (\log (D (h_n f h_n^{-1}) )) = \mathrm{var} (\log ( (Df^n)^{1/n})) = \frac{1}{n} \mathrm{var} (\log (Df^n)),$$
and the right-hand-side expression converges to zero because of the hypothesis of vanishing of the asymptotic distortion.
$\hfill\square$

\vspace{0.3cm}

\begin{rem} The vanishing of the asymptotic distortion of a diffeomorphism that approaches the rotation by conjugates in the $C^{1+bv}$ 
topology can be also seen as a consequence of the (straightforward to establish) fact that the asymptotic distortion is upper semicontinuous 
(as a function defined on the space of $C^{1+bv}$ diffeomorphisms), together with the obvious fact that rotations have vanishing asymptotic distortion.
\end{rem}

\begin{rem} If $f$ is a $C^{1+bv}$ circle diffeomorphism of irrational rotation number and nontrivial image under the Mather homomorphism 
$\mathcal{M}$ (see \S \ref{section-NVAD}), then it cannot approach the rotation under conjugates. This directly follows from the facts that $\mathcal{M}$ 
is conjugacy invariant and continuous. However, it can also be seen as a consequence of Theorem 3 and the fact that nontrivial image under the 
Mather homomorphism implies nonvanishing of the asymptotic distortion (this was established along the proof of Proposition \ref{prop-NVAD}).
\end{rem}

\begin{rem} The preceding proposition has a version for piecewise-affine homeomorphisms. Namely, if $f$ is a such a map with 
vanishing asymptotic distortion, then there exists a sequence $h_n$ of piecewise-affine homeomorphisms such that the total variation
 of the logarithm of the derivative of $h_n f h_n^{-1}$ converges to zero (just follow the same arguments of proof and notice that the 
 explicit formula (\ref{explicit-Z}) yields to piecewise-affine conjugating maps). Remind that, for the case of the circle, the vanishing of 
 the asymptotic distortion is equivalent to that either the map has finite order or, in case of irrational rotation number, that the complete 
 jump along orbits is always 1 (see Remark \ref{rem-mina} and notice that no nontrivial piecewise-affine homeomorphism of the 
interval can have zero asymptotic distortion due to the presence of semi-hyperbolic fixed points, which obviously force 
positive asymptotic distortion).
\end{rem}

\vspace{0.35cm}

We have showed in Proposition \ref{prop-ADV} that the asymptotic distortion of a $C^{1+bv}$ circle diffeomorphism with absolute continuous 
derivative and irrational rotation number is equal to zero. In case of rational rotation number, this is no longer true since  hyperbolic periodic 
points obviously force positive asymptotic distortion. Actually, as we next see, hyperbolicity is the only obstruction to the vanishing of the 
asymptotic distortion for maps in $\mathrm{PSL} (2,\mathbb{R})$, yet it is not the only one in general.

\vspace{0.1cm}

\begin{prop} \label{prop-mobius}
A M\"obius transformation has positive asymptotic distortion if and only if it is hyperbolic.
\end{prop}

\vspace{0.1cm}

\noindent{\bf Proof.} There are (at least) two ways to prove this. First, since hyperbolic M\"obius transformations have hyperbolic fixed points 
and, hence, positive asymptotic distortion, we need to check that both elliptic and parabolic M\"obius transformations have zero asymptotic 
distortion. This can be showed using Theorem 3. Namely, elliptic M\"obius transformations are conjugate to Euclidean rotations, and parabolic 
M\"obius maps admit sequences of conjugates that converge to the identity, as it readily follows by letting $\lambda \to 0$ in the identity
$$ \left[ {\begin{array}{cc}
   \lambda & 0 \\
   0 & 1 / \lambda \\
  \end{array} } \right] 
  \left[ {\begin{array}{cc}
   1 & t \\
   0 & 1 \\
  \end{array} } \right]
  \left[ {\begin{array}{cc}
   \lambda & 0 \\
   0 & 1 / \lambda \\
  \end{array} } \right]^{-1} =  
  \left[ {\begin{array}{cc}
   1 & \lambda^2 t \\
   0 & 1 \\
  \end{array} } \right] .
  $$

The starting observation for the second proof is the following: consider a M\"obius transformation $f \!\in\! \mathrm{PSL}(2,\mathbb{R})$ 
both as an isometry of the Poincar\'e disk (endowed with the hyperbolic distance $\mathrm{dist}_h$) and as a circle diffeomorphism. 
Then the following equality holds: 
\begin{equation}\label{adolfo}
\mathrm{var} (\log (Df); \clo) = 4 \, \mathrm{dist}_h (f (0), 0).
\end{equation}
Assuming this equality, we conclude that the rate of $\, \mathrm{var} (\log (Df^n); \clo) \,$ is linear if and only if the same happens 
for $\, \mathrm{dist}_h (f^n(0),0), \,$ and it is very well known that the behavior of the last sequence distinguishes M\"obius maps: 
this sequence is bounded, unbounded but sublinear, or linear, in the elliptic, parabolic, or hyperbolic case, respectively. 

We are hence left to prove (\ref{adolfo}), which is just a computation. Namely, if we write $f$ in complex form as
$$f (e^{i\theta}) = e^{i\alpha} \cdot \frac{e^{i\theta} - a}{1 - \bar{a} e^{i\theta} }, \qquad a = re^{i\sigma}, \,\, |r| < 1,$$
then, as a circle diffeomorphism, $df / d\theta$ is the norm of the derivative of this expression, which equals
$$
\left\| \frac{e^{i\alpha} [i e^{i\theta} (1 - re^{i(\theta - \sigma)}) + ( e^{i\theta} - re^{i\alpha}) r i e^{i(\theta - \alpha)}]}{(1 - r e^{i (\theta - \sigma)})^2} \right\| 
= 
\frac{1 - r^2}{1 - 2 r \cos (\theta - \sigma) + r^2}.$$
The last expression attains its maximum at $\theta = \sigma + \pi/2$ and its minimum at $\theta = \sigma - \pi/2$, and it is monotone 
(as a function of $\theta$) on both intervals with these endpoints. Using this, one easily concludes that
$$\mathrm{var} (\log (Df); \clo) 
= 2 \left[ \log \left( \frac{1 - r^2}{(1+r)^2} \right) - \log \left( \frac{1-r^2}{(1-r)^2} \right ) \right] 
= 4 \log \left( \frac{1 - r}{1 + r} \right).$$
Finally, it is well known that the last expression corresponds to \, $4 \, \mathrm{dist}_h (0, -e^{i\alpha} a)$, which closes the proof since 
$f (0) = -e^{i\alpha} a$. $\hfill\square$

\vspace{0.5cm}

Despite Proposition \ref{prop-mobius} above, the naive conjecture that, for general diffeomorphisms, positive asymptotic distortion 
implies the existence of hyperbolic periodic points is false. Indeed, below we provide an example of a $C^{\infty}$ 
diffeomorphism of $[0,1]$ having no fixed point in the interior and with positive asymptotic distortion yet the endpoints 
are parabolic fixed points. This slightly extends a claim by Farinelli \cite{farinelli}: there exist $C^{1+vb}$ diffeomorphisms 
of the interval that do not approach the identity by $C^{1+bv}$ conjugates though their fixed points are all parabolic (Farinelli's 
claim concerns $C^2$ regularity, where a similar statement follows as a direct application of Mather's invariant \cite{yoccoz-ast}).

\vspace{0.15cm}

\begin{ex} Start with a parabolic M\"obius transformation, and denote $\hat{f}$ its restriction 
to a fixed interval having no fixed point inside. After conjugacy by an affine map, we may assume that this interval 
is $[0,1]$ and (passing to the inverse if necessary) that $\hat{f}(x) > x$ for all $x \in (0,1)$. By Proposition \ref{prop-mobius} above, the 
asymptotic distortion of $\hat{f}$ vanishes. 

Fix $a < b < \hat{f} (a)$ in $(0,1)$, and let $f$ be a $C^{\infty}$ diffeomorphism with no fixed point in 
the interior that coincides with $\hat{f}$ on a neighborhood 
of $[0,a] \cup [f(a), 1]$ and such that $f(b) = \hat{f} (b)$ but $Df (b) \neq D\hat{f} (b)$.. 
We claim that $f$ has nonzero asymptotic distortion. Indeed, since the asymptotic distortion 
of $\hat{f}$ is zero, given $\varepsilon > 0$ there exists $N$ such that, for all $n \geq N$,
$$\sum_{i = 1}^{n} | \log (D \hat{f}^n (\hat{f}^{-n+i} (b))) - \log (D\hat{f}^n (\hat{f}^{-n+i} (a))) | \leq \mathrm{var} (\log(D \hat{f}^n); [0,1]) 
\leq n \, \varepsilon.$$
Fix \, $\varepsilon < \big| \log (D\hat{f} (b)) - \log (Df(b)) \big| / 2$. \, The preceding inequality may be rewritten as 
$$\sum_{i = 1}^{n} \left| \sum_{j=0}^{n-1} \log (D \hat{f} (\hat{f}^{-n+i+j} (b))) - \log (D\hat{f} (\hat{f}^{-n+i+j} (a))) \right| \leq n \, \varepsilon.$$
In the left-hand side expression, we may replace $\hat{f}$ by $f$ at each occurrence except for that, for each index $i$, the value of 
$Df$ is different from that of $D\hat{f}$ at the point $f^{-n+i+j} (b) = \hat{f}^{-n+i+j} (b)$ for $j = n-i$. This yields the inequality
$$\sum_{i = 1}^{n} \left| [\log(D \hat{f} (b)) - \log(Df(b))] + \sum_{j=0}^{n-1} \log (D f (f^{-n+i+j} (b))) - \log (D f (f^{-n+i+j} (a))) \right| 
\leq n \, \varepsilon.$$
By the triangle inequality $|y| \geq |x| - |x+y|$, this implies that 
$$\sum_{i = 1}^{n} \left| \sum_{j=0}^{n-1} \log (D f (f^{-n+i+j} (b))) - \log (D f (f^{-n+i+j} (a))) \right| 
\geq n \, \big| \log(D \hat{f} (b)) - \log(Df(b)) \big| - n \, \varepsilon,$$
Since the left-hand side expression above may be rewritten as 
$$\sum_{i=1}^n \big| \log(Df^n (f^{-n+i}(b))) - \log(Df^n (f^{-n+i}(a))) \big|,$$
the choice of $\varepsilon$ implies that 
$$\mathrm{var} (\log(Df^n); [0,1]) \geq \frac{n \, | \log(D \hat{f} (b)) - \log(Df(b)) |}{2}.$$
Letting $n$ go to infinite, this shows that 
$$dist_{\infty} (f) \geq \frac{| \log(D \hat{f} (b)) - \log(Df(b)) |}{2}.$$
\end{ex}

\vspace{0.15cm}

\begin{rem} The maps $\hat{f}$ and $f$ cannot be conjugated by a $C^{1+bv}$ diffeomorphism, since 
one of them has vanishing asymptotic distortion and the other one does not. A less obvious claim is that they are 
non $C^1$ conjugate, as we show below. The reader will find in \cite{EN} a more precise result in this direction.

Assume that $h$ is a $C^1$ diffeomorphism satisfying $h \hat{f} h^{-1} = f$. Since $\hat{f} = f$ close to the endpoints, the germs of $h$ 
and $f$ at these endpoints commute. By Kopell's lemma (see \cite{libro}), if we denote by $\hat{f}^t$ the flow associated to $\hat{f}$,  
there exist $t_0,t_1$ such that, close to the origin (resp. close to $1$), one has 
$h = \hat{f}^{t_0}$ (resp. $h = \hat{f}^{t_1}$). If we change $h$ by $h (\hat{f}^{t_0})^{-1}$, we still have the conjugacy relation $\, h \hat{f} h^{-1} = f, \,$ 
but now we have that $h$ equals the identity close to the origin. Since $\hat{f}$ and $f$ coincide on $[0,a] \cup [g (a),1]$ and both 
have the same action on $b$, it is easy to check the following:

\begin{itemize}

\item$h$ equals the identity on $[0,a]$;

\item if we denote by $\bar{h}$ the restriction of $f \hat{f}^{-1}$ to $[a,f(a)]$, then on each interval $f^n ([a,b])$, 
the map $h$ equals $\hat{f}^{n-1} \bar{h} \hat{f}^{-(n-1)}$;

\item $h$ fixes all points of the form $\hat{f}^n (a) = f^n (a)$ and $\hat{f}^n (b) = f^n (b)$.

\end{itemize}

\noindent Now notice that, on the one hand,  
$$Dh (1) = \lim_{n \to \infty} \frac{h(1) - h (f^n(b))}{1 - f^n (b)} = 1.$$
On the other hand, using the chain rule we obtain that, for all $n \geq 1$,
$$D h (f^{n} (b)) = D \bar{h} (b) = \frac{D f (b)}{D \hat{f} (b)} \neq 1.$$
However, these conclusions are incompatible with the continuity of the derivative of $h$ at 1.
\end{rem}

\vspace{0.25cm}

Starting with arbitrary nontrivial diffeomorphisms of zero asymptotic distortion, the technique above easily shows 
that the set of diffeomorphisms with positive asymptotic distortion is dense in the (closed) set of diffeomorphisms 
having only parabolic fixed points. Since the former set coincides with 
$$\bigcup_m \bigcap_n \left\{ f \! : \frac{\mathrm{var}(\log( D f^n); \clo)}{n} \geq \frac{1}{m} \right\},$$
this shows that positive asymptotic distortion is a generic property in the latter set. However,  the construction reveals 
even more than this: zero asymptotic distortion seems to force a very particular structure for the diffeomorphism. 
Despite this, we were not able to settle the question below:

\vspace{0.5cm}

\noindent{\bf Question 3.} Does there exist a real-analytic diffeomorphism of the interval with no hyperbolic 
fixed point but having positive asymptotic distortion?\footnote{This was recently answered in the affirmative by Cohen 
in \cite{cohen}; see also \cite{EN} for a more systematic study.}

\vspace{0.5cm}

Many other natural questions concerning the value of the asymptotic distortion for specific types of diffeomorphisms arise.

\vspace{0.5cm}

\noindent{\bf Question 4.} Can the asymptotic distortion of a nontrivial piecewise-projective 
homeomorphism that is not piecewise projectively conjugate to a parabolic M\"obius map vanish?

\vspace{0.5cm}

Remind that the asymptotic distortion of every nontrivial piecewise-affine homeomorphism of the interval is positive, 
due to the presence of semi-hyperbolic fixed points.

\vspace{0.5cm}

\noindent{\bf Question 5.} What is the range of values of the asymptotic distortion for the $C^{\infty}$ diffeomorphisms 
of the interval that arise in the Ghys-Sergiescu's smooth realizations of Thompson's group F (see \cite{GS})? Is it 
possible in this framework to produce (nontrivial) diffeomorphisms with the same asymptotic distortion which 
are not smoothly conjugate?

\vspace{0.5cm}

Another natural question (inspired by \cite{Polt-Sod}) deals with the behavior of the sequence of 
distortions along iterates. (Compare Question 1.)

\vspace{0.5cm}

\noindent{\bf Question 6.} Assume that the asymptotic distortion of a diffeomorphism $f$ of $[0,1]$ with no hyperbolic fixed point 
vanishes. Can the sequence \, $\mathrm{var} (\log (Df^n); [0,1])$ \, grow faster than logarithmically?

\vspace{0.5cm}

Finally, having established the existence of diffeomorphisms with positive asymptotic distortion (besides those with hyperbolic 
periodic points), the next generalization of Theorem~3 becomes interesting. The proof is a 
straightforward modification that we leave to the reader. 

\vspace{0.15cm}

\begin{prop} Let $f$ be a $C^{1+bv}$ diffeomorphism of a compact manifold $M$. If there is a sequence 
$h_n$ of $C^{1+bv}$ diffeomorphisms of $M$ such that, for some $d \geq 0$ and all integers $n \geq 1$, one has
$$\mathrm{var} ( \log (D (h_n f h_n^{-1})); M) < d + \frac{1}{n},$$
then \, $d_{\infty}(f) \leq d$. \, Conversely, for each $n \geq 1$ 
there exists a $C^{1+bv}$ diffeomorphism of $M$ such that
$$\mathrm{var} (\log (D (h_n f h_n^{-1})); M) < d_{\infty} (f) + \frac{1}{n}.$$
\end{prop}

\vspace{0.5cm}

In other words, $\mathrm{dist}_{\infty} (f)$ is nothing but the infimum of the variation of the logarithm of the 
derivatives of conjugates of $f$, that is, a dynamically invariant version of $\mathrm{var} (\log (D (\cdot)))$. 

We close with a question relating asymptotic distortion and centralizers of $C^{1+bv}$ diffeomorphisms 
of the interval.\footnote{This was recently answered in the negative in \cite{EN}.}

\vspace{0.5cm}

\noindent{\bf Question 7.} Can a nontrivial $C^{\infty}$ diffeomorphism $f$ of $[0,1]$ with no fixed point in the interior and 
positive asymptotic distortion be centralized by a $C^{\infty}$ diffeomorphism $g$ of zero asymptotic distortion in a way 
that $\langle f,g \rangle \sim \mathbb{Z}^2$~?


\section{On $C^2$ conjugates of commuting diffeomorphisms}
\label{section-C2}

Our Main Theorem gives no information on the $C^2$ closure of the conjugacy class of $C^2$ circle diffeomorphisms. 
Whether or not they contain the corresponding rotation in case of irrational rotation number seems to be a difficult question. 
Theorem 4 (recalled below) is a partial result in this direction which is worth to include, despite the argument of proof is rather standard. 

\vspace{0.5cm}

\noindent{\bf Theorem 4.} {\em If $\Gamma$ is a finitely-generated Abelian group of $C^2$ circle diffeomorphisms that is conjugate 
to a group of rotations by a $C^1$ diffeomorphism, then there exists a sequence $h_n$ of $C^2$ diffeomorphisms such that, for all 
$f \in \Gamma$, the sequence of conjugates $g_nfg_n^{-1}$ converges to a rotation in $C^2$ topology as $n$ goes to infinity.}

\vspace{0.2cm}

\noindent{\bf Proof.} We start by reducing the general case to that of torsion-free groups. To do this, notice that torsion elements of $\Gamma$ 
form a finite cyclic subgroup, say generated by an element  $f$ of order $n$. If we conjugate $\Gamma$ by the diffeomorphism $h$ defined by
$$h := \frac{1}{n} \left[ Id + f + \ldots + f^{n-1} \right],$$
then this finite group becomes a group $\Gamma_*$ of Euclidean rotations. The quotient of $\mathrm{S}^1$ 
under the action of $\Gamma_*$ is a smooth circle endowed with a $\Gamma$-action in which torsion elements act trivially. After proving 
the claim of the theorem for this action, a lift of the conjugating maps to $\mathrm{S^1}$ will commute with the elements of $\Gamma_*$, 
thus providing the desired conjugacies. 

Assume henceforth that $\Gamma$ is isomorphic to $\mathbb{Z}^k$. Let $\{ f_1, \ldots, f_k\}$ be a generating system of $\Gamma$, and for each 
$n$ denote $B(n)$ the set of elements of the form $f_1^{n_1} f_2^{n_2} \cdots f_k^{n_k}$ for which $0 \leq n_i \leq n$. As in the Main Proposition, 
for each $n \geq 1$, we let $h_n$ be a diffeomorphism satisfying (\ref{explicit-general}). For each $\theta$, letting $\sigma := \varphi_n^{-1} (\theta)$, 
a direct computation yields, for each $1 \leq i \leq k$,
\begin{eqnarray*}
D (h_n \circ f_i \circ h_n^{-1}) (\theta)
&=&
D h_n (f_i (\sigma)) \cdot D f_i (\sigma) \cdot D \varphi_n^{-1} (\theta) \\
&=&
D f_i (\sigma) \cdot \frac{D h_n (f_i (\sigma)))}{D \varphi_n (\sigma)} \\
&=&
D f_i (\sigma) \cdot \left[ \frac{\prod_{f \in B(n-1)} Df (f_i (\sigma)) }{\prod_{f \in B(n-1)} Df (\sigma)} \right]^{1 / |B(n-1)|} \\
&=&
D f_i (\sigma) \cdot \left[ \frac{\prod_{0 \leq n_j < n} D (f_1^{n_1} \cdots f_k^{n_k}) (f_i (\sigma)) }{\prod_{0 \leq n_j < n} D (f_1^{n_1} \cdots f_k^{n_k}) (\sigma)} \right]^{1 / |B(n-1)|} \\
&=&
\left[ \frac{\prod_{0 \leq n_j < n} D (f_1^{n_1} \cdots f_k^{n_k}) (f_i (\sigma)) \cdot Df_i (\sigma)}{\prod_{0 \leq n_j < n} D (f_1^{n_1} \cdots f_k^{n_k}) (\sigma)} \right]^{1 / |B(n-1)|} \\
\end{eqnarray*}
Using the chain rule, a telescopic multiplication yields
$$D (h_n \circ f_i \circ h_n^{-1}) (\theta)
= 
\left[ 
\frac{\prod_{0 \leq n_j < n, j \neq i} D (f_1^{n_1} \cdots f_{i-1}^{n_{i-1}} f_i^n f_{i+1}^{n_{i+1}} \cdots f_k^{n_k}) (\sigma) }
{\prod_{0 \leq n_j < n, j \neq i} D (f_1^{n_1} \cdots f_{i-1}^{n_{i-1}}  f_{i+1}^{n_{i+1}} \cdots f_k^{n_k}) (\sigma)} 
\right]^{1 / |B(n-1)|}.$$
Using the identity $D^2 (g) / D (g) = D (\log (D(g)))$, we obtain that 
$$\frac{D^2 (h_n \circ f_i \circ h_n^{-1}) (\theta)}{D (h_n \circ f_i \circ h_n^{-1}) (\theta)}$$
equals
\begin{small}
$$\frac{1}{|B (n-1)|} \sum_{\substack{0 \leq n_j < n\\ j \neq i}} \! \left[ \frac{D^2 (f_1^{n_1} \cdots f_i^n \cdots f_k^{n_k}) (\sigma)}
{D (f_1^{n_1} \cdots  f_i^n \cdots f_k^{n_k}) (\sigma)} \cdot D h_n^{-1} (\theta) - 
\frac{D^2 (f_1^{n_1} \cdots f_{i-1}^{n_{i-1}} f_{i+1}^{n_{i+1}} \cdots f_k^{n_k}) (\sigma)}
{D (f_1^{n_1} \cdots f_{i-1}^{n_{i-1}} f_{i+1}^{n_{i+1}} \cdots f_k^{n_k}) (\sigma)} \cdot D h_n^{-1} (\theta) \right] \!.$$
\end{small}The cocycle relation (\ref{cociclo}) 
allows transforming the last expression into 
$$\frac{D h_n^{-1} (\theta)}{|B (n\!-\!1)|} \sum_{\substack{0 \leq n_j < n\\ j \neq i}} 
\frac{D^2 f_i^n (f_1^{n_1} \cdots f_{i-1}^{n_{i-1}} f_{i+1}^{n_{i+1}} \cdots f_k^{n_k} (\sigma))}
{D f_i^n (f_1^{n_1} \cdots f_{i-1}^{n_{i-1}} f_{i+1}^{n_{i+1}} \cdots f_k^{n_k} (\sigma))} 
\cdot 
D (f_1^{n_1} \cdots f_{i-1}^{n_{i-1}} f_{i+1}^{n_{i+1}} \cdots f_k^{n_k}) (\sigma),$$
which, in its turn, transforms into 
\begin{small}
\begin{equation}\label{a-controlar}
\frac{D h_n^{-1} (\theta)}{|B (n\!-\!1)|} \sum_{\substack{0 \leq n_j < n\\ j \neq i}} \sum_{k=0}^{n-1}
\frac{D^2 f_i (f_i^k f_1^{n_1} \cdots f_{i-1}^{n_{i-1}} f_{i+1}^{n_{i+1}} \cdots f_k^{n_k} (\sigma))}
{D f_i (f_i^k f_1^{n_1} \cdots f_{i-1}^{n_{i-1}} f_{i+1}^{n_{i+1}} \cdots f_k^{n_k} (\sigma))} 
\cdot
D (f_i^k f_1^{n_1} \cdots f_{i-1}^{n_{i-1}} f_{i+1}^{n_{i+1}} \cdots f_k^{n_k}) (\sigma).
\end{equation}
\end{small}

Recall that we are assuming that $\Gamma$ is $C^1$ conjugate to a group of rotations. 
Let hence $\varphi$ be a $C^1$ diffeomorphism  such that, for all $f \in \Gamma$,
\begin{equation}\label{conja}
f = \varphi \circ R_{\rho} \circ \varphi^{-1}, \qquad \rho = \rho(f).
\end{equation}
The derivative 
\begin{equation}\label{de-efe}
D f (\nu) = \frac{D \varphi (\varphi^{-1}(\nu) + \rho)}{D \varphi (\varphi^{-1} (\nu))}
\end{equation}
is uniformly bounded by $C := \max D \varphi / \min D \varphi$ and, by definition (\ref{explicit-general}), the same holds for the derivative 
of $h_n$. Passing to the inverse, we obtain the same upper bound for $D h_n^{-1}$, which allows to control the 
numerator of the left factor of (\ref{a-controlar}). Concerning the sum in (\ref{a-controlar}), notice that it is made by sums of the form
$$\sum_{k=0}^{n-1}
\frac{D^2 f_i (f_i^k f_1^{n_1} \cdots f_{i-1}^{n_{i-1}} f_{i+1}^{n_{i+1}} \cdots f_k^{n_k} (\sigma))}
{D f_i (f_i^k f_1^{n_1} \cdots f_{i-1}^{n_{i-1}} f_{i+1}^{n_{i+1}} \cdots f_k^{n_k} (\sigma))} 
\cdot
D (f_i^k f_1^{n_1} \cdots f_{i-1}^{n_{i-1}} f_{i+1}^{n_{i+1}} \cdots f_k^{n_k}) (\sigma)$$
in a number that equals the cardinal of sets of integers $n_j$ satisfying $0 \leq n_j < n, j \neq i$. 
Since this cardinal is obviously equal to $|B(n-2)|$ and 
$$ \frac{ |B(n-2)| }{ |B(n-1)| } = \frac{1}{n},$$ 
the proof will be finished if we show that the expressions 
$$\frac{1}{n} \sum_{k=0}^{n-1}
\frac{D^2 f_i (f_i^k f_1^{n_1} \cdots f_{i-1}^{n_{i-1}} f_{i+1}^{n_{i+1}} \cdots f_k^{n_k} (\sigma))}
{D f_i (f_i^k f_1^{n_1} \cdots f_{i-1}^{n_{i-1}} f_{i+1}^{n_{i+1}} \cdots f_k^{n_k} (\sigma))} 
\cdot
D (f_i^k f_1^{n_1} \cdots f_{i-1}^{n_{i-1}} f_{i+1}^{n_{i+1}} \cdots f_k^{n_k}) (\sigma)$$ 
uniformly converge to zero as $n \to \infty$. Now, denoting $\rho_i$ the rotation number of $f_i$ and letting $\eta := \varphi^{-1} (\sigma)$, 
relations (\ref{conja}) and (\ref{de-efe}) allow to transform the last mean into
$$\frac{1}{n} \sum_{k=0}^{n-1} \frac{D^2 f_i (\varphi (\eta + k \rho_i + \sum_{j \neq i} n_j \rho_j))}{D f_i (\varphi (\eta + k \rho_i + \sum_{j \neq i} n_j \rho_j))} 
\cdot \frac{D \varphi (\eta + k\rho_i + \sum_{j \neq i} n_j \rho_j) } {D \varphi (\sigma)}.$$
By Weil's equidistribution theorem, this converges as $n \to \infty$ to 
$$\frac{1}{D \varphi (\sigma)} \int_{\mathrm{S}^1} \frac{D^2 f_i (\varphi(\eta))}{D f_i (\varphi(\eta))} \cdot D \varphi (\eta) \, d \eta 
= 
\frac{1}{D \varphi (\sigma)} \int_{\mathrm{S}^1} \frac{D^2 f_i (\sigma)}{D f_i (\sigma)} \, d \sigma = 0,$$
which closes the proof.
$\hfill\square$

\vspace{0.25cm}

\begin{rem}
As in the end of the proof of the Main Theorem, it is not hard to show existence of paths of conjugate actions rather than sequences of 
conjugate actions converging to actions by rotations. Moreover, an exhaustion type argument allows extending Theorem 4 to countable 
Abelian groups. Details are left to the reader.
\end{rem}

\vspace{0.15cm}

\noindent{\bf Acknowledgments.} I'm strongly indebted to Adolfo Guillot, with whom formula (\ref{adolfo}) arose in a discussion during a 
seminar (more than ten years ago). I'm also indebted to Jairo Bochi, Bassam Fayad, Nancy Guelman,  Alejandro Kocsard and Rapha\"el 
Krikorian for quite stimulating discussions, as well as to the two referees whose remarks allowed greatly improving the exposition. 

The preparation of this text was funded by the CONICYT Projects MathAmsud 17-03 
``GDAR: geometry, dynamics, and Anosov representations" and FONDECYTs 1160541 and 1200114 (in Chile), as well as  the  
Project FORDECYT 265667 and the PREI Program of the DGAPA at UNAM (in M\'exico).

\vspace{0.2cm}


\begin{footnotesize}

\vspace{0.1cm}

\noindent Andr\'es Navas\\ 

\vspace{0.2cm}

\noindent Dpto de Matem\'atica y C.C., Univ. de Santiago de Chile\\ 

\noindent Alameda 3363, Estaci\'on Central, Santiago, Chile\\ 

\vspace{0.2cm}

\noindent Unidad Cuernavaca Instituto de Matem\'aticas\\

\noindent Universidad Nacional Aut\'onoma de M\'exico, Campus Morelos\\

\vspace{0.2cm}

\noindent Email address: andres.navas@usach.cl \\ 

\end{footnotesize}

\end{document}